\documentclass{amsart}
\usepackage{amsmath}
\usepackage{amssymb}
\usepackage{amsfonts}

\newtheorem{theorem}{Theorem}
\theoremstyle{plain}

\newtheorem{corollary}{Corollary}

\newtheorem{definition}{Definition}

\newtheorem{lemma}{Lemma}

\numberwithin{equation}{section}
\input{tcilatex}

\begin{document}
\title[Ostrowski type inequalities]{Inequalities of Ostrowski's type for $m-$
and $(\alpha ,m)-$ logarithmically convex functions via Riemann-Liouville
Fractional integrals}
\author{Ahmet Ocak Akdemir}
\address{A\u{g}r\i\ \.{I}brahim \c{C}e\c{c}en University, Faculty of Science
and Letters, Department of Mathematics, A\u{g}r\i , 04100, Turkey.}
\email{ahmetakdemir@agri.edu.tr}

\begin{abstract}
In this paper, we establish some new Ostrowski's type inequalities for $m-$
and $(\alpha ,m)-$ logarithmically convex functions by using the
Riemann-Liouville fractional integrals.
\end{abstract}

\maketitle

\section{INTRODUCTION}

Let $f:I\subset \left[ 0,\infty \right] \rightarrow 
\mathbb{R}
$ be a differentiable mapping on $I^{\circ }$, the interior of the interval $%
I$, such that $f^{\prime }\in L\left[ a,b\right] $ where $a,b\in I$ with $%
a<b $. If $\left\vert f^{\prime }\left( x\right) \right\vert \leq M$, then
the following inequality holds (see \cite{10}).

\begin{equation}
\left\vert f(x)-\frac{1}{b-a}\int_{a}^{b}f(u)du\right\vert \leq \frac{M}{b-a}%
\left[ \frac{\left( x-a\right) ^{2}+\left( b-x\right) ^{2}}{2}\right]
\label{h.1.1}
\end{equation}

This inequality is well known in the literature as the Ostrowski inequality%
\textit{.}\textbf{\ }For some results which generalize, improve and extend
the inequality (\ref{h.1.1}) see (\cite{10, a1, 2, 11, a2, a}) and the
references therein.

Let us recall some known definitions and results which we will use in this
paper. The function $f:\left[ a,b\right] \rightarrow 
\mathbb{R}
,$ is said to be convex, if we have%
\begin{equation*}
f\left( tx+\left( 1-t\right) y\right) \leq tf\left( x\right) +\left(
1-t\right) f\left( y\right)
\end{equation*}%
for all $x,y\in \left[ a,b\right] $ and $t\in \left[ 0,1\right] .$ We can
define starshaped functions on $\left[ 0,b\right] $ which satisfy the
condition%
\begin{equation*}
f\left( tx\right) \leq tf\left( x\right)
\end{equation*}%
for $t\in \left[ 0,1\right] .$

The concept of $m-$convexity has been introduced by Toader in \cite{TOA}, an
intermediate between the ordinary convexity and starshaped property, as
following:

\begin{definition}
The function $f:\left[ 0,b\right] \rightarrow 
\mathbb{R}
,$ $b>0,$ is said to be $m-$convex, where $m\in \left[ 0,1\right] ,$ if we
have%
\begin{equation*}
f\left( tx+m\left( 1-t\right) y\right) \leq tf\left( x\right) +m\left(
1-t\right) f\left( y\right)
\end{equation*}%
for all $x,y\in \left[ 0,b\right] $ and $t\in \left[ 0,1\right] .$ We say
that $f$ is $m-$concave if $-f$ is $m-$convex.
\end{definition}

In \cite{MIH}, Mihe\c{s}an gave definition of $(\alpha ,m)-$convexity as
following;

\begin{definition}
The function $f:[0,b]\rightarrow 
\mathbb{R}
,$ $b>0$ is said to be $(\alpha ,m)-$convex, where $(\alpha ,m)\in \lbrack
0,1]^{2},$ if we have%
\begin{equation*}
f(tx+m(1-t)y)\leq t^{\alpha }f(x)+m(1-t^{\alpha })f(y)
\end{equation*}%
for all $x,y\in \lbrack 0,b]$ and $t\in \lbrack 0,1].$
\end{definition}

Denote by $K_{m}^{\alpha }(b)$ the class of all $(\alpha ,m)-$convex
functions on $[0,b]$ for which $f(0)\leq 0.$ If we choose $(\alpha ,m)=(1,m)$%
, it can be easily seen that $(\alpha ,m)-$convexity reduces to $m-$%
convexity and for $(\alpha ,m)=(1,1),$ we have ordinary convex functions on $%
[0,b].$ For the recent results based on the above definitions see the papers 
\cite{BPR}-\cite{2}.

\begin{definition}
(\cite{1}) A function $f:[0,b]\rightarrow (0,\infty )$ is said to be $m-$%
logarithmically convex if the inequality%
\begin{equation}
f\left( tx+m\left( 1-t\right) y\right) \leq \left[ f\left( x\right) \right]
^{t}\left[ f\left( y\right) \right] ^{m\left( 1-t\right) }  \label{11}
\end{equation}%
holds for all $x,y\in \lbrack 0,b]$, $m\in (0,1]$, and $t\in \lbrack 0,1]$.
\end{definition}

Obviously, if putting $m=1$ in Definition 3, then $f$ is just the ordinary
logarithmically convex function on $\left[ 0,b\right] $.

\begin{definition}
(\cite{1}) A function $f:[0,b]\rightarrow (0,\infty )$ is said to be $\left(
\alpha ,m\right) -$logarithmically convex if%
\begin{equation}
f\left( tx+m\left( 1-t\right) y\right) \leq \left[ f\left( x\right) \right]
^{t^{\alpha }}\left[ f\left( y\right) \right] ^{m\left( 1-t^{\alpha }\right)
}  \label{2}
\end{equation}%
holds for all $x,y\in \lbrack 0,b]$, $\left( \alpha ,m\right) \in \left( 0,1%
\right] \times \left( 0,1\right] ,$ and $t\in \lbrack 0,1]$.
\end{definition}

Clearly, when taking $\alpha =1$ in Definition 4, then $f$ becomes the
standard $m$-logarithmically convex function on $\left[ 0,b\right] $.{}

We give some necessary definitions and mathematical preliminaries of
fractional calculus theory which are used throughout this paper.

\begin{definition}
Let $f\in L_{1}[a,b].$ The Riemann-Liouville integrals $J_{a^{+}}^{\mu }f$
and $J_{b^{-}}^{\mu }f$ of order $\mu >0$ with $a\geq 0$ are defined by%
\begin{equation*}
J_{a^{+}}^{\mu }f\left( x\right) =\frac{1}{\Gamma (\mu )}\underset{a}{%
\overset{x}{\int }}\left( x-t\right) ^{\mu -1}f(t)dt,\text{ \ }x>a
\end{equation*}%
and%
\begin{equation*}
J_{b^{-}}^{\mu }f\left( x\right) =\frac{1}{\Gamma (\mu )}\underset{x}{%
\overset{b}{\int }}\left( t-x\right) ^{\mu -1}f(t)dt,\text{ \ }x<b
\end{equation*}%
respectively where $\Gamma (\mu )=\underset{0}{\overset{\infty }{\int }}%
e^{-t}u^{\mu -1}du.$ Here is $J_{a^{+}}^{0}f(x)=J_{b^{-}}^{0}f(x)=f(x).$
\end{definition}

In the case of $\mu =1$, the fractional integral reduces to the classical
integral. For some recent results connected with \ fractional integral
inequalities see \cite{a2}-\cite{a}.

The aim of this study is to establish some Ostrowski type inequalities for
the class of functions whose derivatives in absolute value are $m-$ and $%
\left( \alpha ,m\right) -$ geometrically convex functions via
Riemann-Liouville fractional integral.

\section{THE\ NEW\ RESULTS}

In order to prove our results, we need the following lemma that has been
obtained in \cite{a2}:

\begin{lemma}
(\cite{a2})\label{L1} Let $f:\left[ a,b\right] \rightarrow 
\mathbb{R}
$ be a differentiable mapping on $\left( a,b\right) $ with $a<b.$ If $%
f^{\prime }\in L\left[ a,b\right] ,$ then for all $x\in \left[ a,b\right] $
and $\mu >0$ we have:%
\begin{eqnarray*}
&&\frac{\left( x-a\right) ^{\mu }+\left( b-x\right) ^{\mu }}{b-a}f\left(
x\right) -\frac{\Gamma \left( \mu +1\right) }{b-a}\left[ J_{x^{-}}^{\mu
}f\left( a\right) +J_{x^{+}}^{\mu }f\left( b\right) \right] \\
&=&\frac{\left( x-a\right) ^{\mu +1}}{b-a}\int_{0}^{1}t^{\mu }f^{\prime
}\left( tx+\left( 1-t\right) a\right) dt+\frac{\left( b-x\right) ^{\mu +1}}{%
b-a}\int_{0}^{1}t^{\mu }f^{\prime }\left( tx+\left( 1-t\right) b\right) dt
\end{eqnarray*}%
where $\Gamma (\mu )=\underset{0}{\overset{\infty }{\int }}e^{-t}u^{\mu
-1}du.$
\end{lemma}

\begin{theorem}
Let $f:\left[ 0,\infty \right) \rightarrow \left( 0,\infty \right) $ be \
differentiable mapping with $a,b\in \left[ 0,\infty \right) $ such that $%
a<b. $ If $\left\vert f^{\prime }\left( x\right) \right\vert $ is $\left(
\alpha ,m\right) -$logarithmically convex function with $\left\vert
f^{\prime }\left( x\right) \right\vert \leq M,f^{\prime }\in L\left[ a,b%
\right] ,$ $\left( \alpha ,m\right) \in \left( 0,1\right] \times \left( 0,1%
\right] $ and $\mu >0,$ then the following inequality for fractional
integrals holds:%
\begin{eqnarray}
&&\left\vert \frac{\left( x-a\right) ^{\mu }+\left( b-x\right) ^{\mu }}{b-a}%
f\left( x\right) -\frac{\Gamma \left( \mu +1\right) }{b-a}\left[
J_{x^{-}}^{\mu }f\left( a\right) +J_{x^{+}}^{\mu }f\left( b\right) \right]
\right\vert  \notag \\
&\leq &\left[ \frac{1}{2\mu +1}+K_{1}(\alpha ,m,t)\right] \left[ \frac{%
\left( x-a\right) ^{\mu +1}+\left( b-x\right) ^{\mu +1}}{2\left( b-a\right) }%
\right]  \label{a1}
\end{eqnarray}%
where%
\begin{equation*}
K_{1}(\alpha ,m,t)=\left\{ 
\begin{array}{cc}
\frac{M^{2m}\left( M^{2\alpha -2\alpha m}-1\right) }{\left( 2\alpha -2\alpha
m\right) \ln M} & ,M<1 \\ 
&  \\ 
1 & ,M=1%
\end{array}%
\right. .
\end{equation*}
\end{theorem}

\begin{proof}
By Lemma \ref{L1} and since $\left\vert f^{\prime }\right\vert $ is $\left(
\alpha ,m\right) -$logarithmically convex$,$ we can write%
\begin{eqnarray*}
&&\left\vert \frac{\left( x-a\right) ^{\mu }+\left( b-x\right) ^{\mu }}{b-a}%
f\left( x\right) -\frac{\Gamma \left( \mu +1\right) }{b-a}\left[
J_{x^{-}}^{\mu }f\left( a\right) +J_{x^{+}}^{\mu }f\left( b\right) \right]
\right\vert \\
&\leq &\frac{\left( x-a\right) ^{\mu +1}}{b-a}\int_{0}^{1}t^{\mu }\left\vert
f^{\prime }\left( tx+\left( 1-t\right) a\right) \right\vert dt+\frac{\left(
b-x\right) ^{\mu +1}}{b-a}\int_{0}^{1}t^{\mu }\left\vert f^{\prime }\left(
tx+\left( 1-t\right) b\right) \right\vert dt \\
&\leq &\frac{\left( x-a\right) ^{\mu +1}}{b-a}\int_{0}^{1}t^{\mu }\left\vert
f^{\prime }\left( x\right) \right\vert ^{t^{\alpha }}\left\vert f^{\prime
}\left( \frac{a}{m}\right) \right\vert ^{m\left( 1-t^{\alpha }\right) }dt+%
\frac{\left( b-x\right) ^{\mu +1}}{b-a}\int_{0}^{1}t^{\mu }\left\vert
f^{\prime }\left( x\right) \right\vert ^{t^{\alpha }}\left\vert f^{\prime
}\left( \frac{b}{m}\right) \right\vert ^{m\left( 1-t^{\alpha }\right) }dt \\
&\leq &\frac{\left( x-a\right) ^{\mu +1}}{b-a}\int_{0}^{1}t^{\mu
}M^{m+t^{\alpha }\left( 1-m\right) }dt+\frac{\left( b-x\right) ^{\mu +1}}{b-a%
}\int_{0}^{1}t^{\mu }M^{m+t^{\alpha }\left( 1-m\right) }dt.
\end{eqnarray*}%
By using the elemantery inequality $cd\leq \frac{c^{2}+d^{2}}{2},$ we have%
\begin{eqnarray}
&&  \label{k1} \\
&&\left\vert \frac{\left( x-a\right) ^{\mu }+\left( b-x\right) ^{\mu }}{b-a}%
f\left( x\right) -\frac{\Gamma \left( \mu +1\right) }{b-a}\left[
J_{x^{-}}^{\mu }f\left( a\right) +J_{x^{+}}^{\mu }f\left( b\right) \right]
\right\vert  \notag \\
&\leq &\frac{\left( x-a\right) ^{\mu +1}}{b-a}\int_{0}^{1}\frac{t^{2\mu
}+M^{2(m+t^{\alpha }\left( 1-m\right) )}}{2}dt+\frac{\left( b-x\right) ^{\mu
+1}}{b-a}\int_{0}^{1}\frac{t^{2\mu }+M^{2(m+t^{\alpha }\left( 1-m\right) )}}{%
2}dt  \notag \\
&=&\left[ \frac{1}{2\mu +1}+\int_{0}^{1}M^{2(m+t^{\alpha }\left( 1-m\right)
)}dt\right] \left[ \frac{\left( x-a\right) ^{\mu +1}+\left( b-x\right) ^{\mu
+1}}{2\left( b-a\right) }\right] .  \notag
\end{eqnarray}%
If we choose $M=1,$ then 
\begin{equation*}
\int_{0}^{1}M^{2(m+t^{\alpha }\left( 1-m\right) )}dt=1.
\end{equation*}%
If $M<1,$ then $M^{2(m+t^{\alpha }\left( 1-m\right) )}\leq M^{2(m+\alpha
t\left( 1-m\right) )}$, thus%
\begin{equation*}
\int_{0}^{1}M^{2(m+\alpha t\left( 1-m\right) )}dt=\frac{M^{2m}\left(
M^{2\alpha -2\alpha m}-1\right) }{\left( 2\alpha -2\alpha m\right) \ln M}
\end{equation*}%
Which completes the proof.
\end{proof}

\begin{corollary}
Let $f:\left[ 0,\infty \right) \rightarrow \left( 0,\infty \right) $ be \
differentiable mapping with $a,b\in \left[ 0,\infty \right) $ such that $%
a<b. $ If $\left\vert f^{\prime }\left( x\right) \right\vert $ is $m-$%
logarithmically convex function with $\left\vert f^{\prime }\left( x\right)
\right\vert \leq M,f^{\prime }\in L\left[ a,b\right] ,$ $m\in \left( 0,1%
\right] $ and $\mu >0,$ then the following inequality for fractional
integrals holds:%
\begin{eqnarray}
&&\left\vert \frac{\left( x-a\right) ^{\mu }+\left( b-x\right) ^{\mu }}{b-a}%
f\left( x\right) -\frac{\Gamma \left( \mu +1\right) }{b-a}\left[
J_{x^{-}}^{\mu }f\left( a\right) +J_{x^{+}}^{\mu }f\left( b\right) \right]
\right\vert  \notag \\
&\leq &\left[ \frac{1}{2\mu +1}+\frac{M^{2}-M^{2m}}{2\ln M-2m\ln M}\right] %
\left[ \frac{\left( x-a\right) ^{\mu +1}+\left( b-x\right) ^{\mu +1}}{%
2\left( b-a\right) }\right] .  \label{a3}
\end{eqnarray}
\end{corollary}

\begin{proof}
If we take $\alpha =1$ in (\ref{a1}), we get the required result.
\end{proof}

\begin{corollary}
Let $f:\left[ 0,\infty \right) \rightarrow \left( 0,\infty \right) $ be \
differentiable mapping with $a,b\in \left[ 0,\infty \right) $ such that $%
a<b. $ If $\left\vert f^{\prime }\left( x\right) \right\vert $ is
logarithmically convex function with $\left\vert f^{\prime }\left( x\right)
\right\vert \leq M,f^{\prime }\in L\left[ a,b\right] $ and $\mu >0,$ then
the following inequality for fractional integrals holds:%
\begin{eqnarray}
&&\left\vert \frac{\left( x-a\right) ^{\mu }+\left( b-x\right) ^{\mu }}{b-a}%
f\left( x\right) -\frac{\Gamma \left( \mu +1\right) }{b-a}\left[
J_{x^{-}}^{\mu }f\left( a\right) +J_{x^{+}}^{\mu }f\left( b\right) \right]
\right\vert  \notag \\
&\leq &\left[ \frac{1}{2\mu +1}+M^{2}\right] \left[ \frac{\left( x-a\right)
^{\mu +1}+\left( b-x\right) ^{\mu +1}}{2\left( b-a\right) }\right] .
\label{a4}
\end{eqnarray}
\end{corollary}

\begin{proof}
If we take $\alpha =m=1$ in (\ref{k1}), we get the required result.
\end{proof}

\begin{corollary}
Let $f:\left[ 0,\infty \right) \rightarrow \left( 0,\infty \right) $ be \
differentiable mapping with $a,b\in \left[ 0,\infty \right) $ such that $%
a<b. $ If $\left\vert f^{\prime }\left( x\right) \right\vert $ is
logarithmically convex function with $\left\vert f^{\prime }\left( x\right)
\right\vert \leq M$ and $f^{\prime }\in L\left[ a,b\right] ,$ then the
following inequality holds:%
\begin{equation*}
\left\vert f\left( x\right) -\frac{1}{b-a}\int_{a}^{b}f(u)du\right\vert \leq %
\left[ \frac{1}{3}+M^{2}\right] \left[ \frac{\left( x-a\right) ^{2}+\left(
b-x\right) ^{2}}{2\left( b-a\right) }\right] .
\end{equation*}
\end{corollary}

\begin{proof}
If we choose $\mu =1$ in (\ref{a4}), we get the required result.
\end{proof}

\begin{theorem}
Let $f:\left[ 0,\infty \right) \rightarrow \left( 0,\infty \right) $ be \
differentiable mapping with $a,b\in \left[ 0,\infty \right) $ such that $%
a<b. $ If $\left\vert f^{\prime }\left( x\right) \right\vert ^{q}$ is $%
\left( \alpha ,m\right) -$logarithmically convex function with $\left\vert
f^{\prime }\left( x\right) \right\vert ^{q}\leq M,f^{\prime }\in L\left[ a,b%
\right] ,$ $\left( \alpha ,m\right) \in \left( 0,1\right] \times \left( 0,1%
\right] $ and $\mu >0,$ then the following inequality for fractional
integrals holds:%
\begin{eqnarray}
&&  \label{a5} \\
&&\left\vert \frac{\left( x-a\right) ^{\mu }+\left( b-x\right) ^{\mu }}{b-a}%
f\left( x\right) -\frac{\Gamma \left( \mu +1\right) }{b-a}\left[
J_{x^{-}}^{\mu }f\left( a\right) +J_{x^{+}}^{\mu }f\left( b\right) \right]
\right\vert  \notag \\
&\leq &\left( \frac{q-1}{\mu \left( q-p\right) +q-1}\right) ^{\frac{q-1}{q}%
}\left( K_{2}(\alpha ,m,t)\right) ^{\frac{1}{q}}\left[ \frac{\left(
x-a\right) ^{\mu +1}+\left( b-x\right) ^{\mu +1}}{\left( b-a\right) }\right]
\notag
\end{eqnarray}%
where $q>1,$ $0\leq p\leq q$ and 
\begin{equation*}
K_{2}(\alpha ,m,t)=\left\{ 
\begin{array}{cc}
\frac{M^{m}\left( \Gamma \left( \mu p+1\right) -\Gamma \left( \mu p+1,\ln
M^{\alpha \left( m-1\right) }\right) \right) }{\left( \ln M^{\alpha \left(
m-1\right) }\right) ^{\mu p+1}} & ,M<1 \\ 
&  \\ 
\frac{1}{\mu p+1} & ,M=1%
\end{array}%
\right. .
\end{equation*}
\end{theorem}

\begin{proof}
From Lemma 1 and by using the properties of modulus, we have%
\begin{eqnarray*}
&&\left\vert \frac{\left( x-a\right) ^{\mu }+\left( b-x\right) ^{\mu }}{b-a}%
f\left( x\right) -\frac{\Gamma \left( \mu +1\right) }{b-a}\left[
J_{x^{-}}^{\mu }f\left( a\right) +J_{x^{+}}^{\mu }f\left( b\right) \right]
\right\vert \\
&\leq &\frac{\left( x-a\right) ^{\mu +1}}{b-a}\int_{0}^{1}t^{\mu }\left\vert
f^{\prime }\left( tx+\left( 1-t\right) a\right) \right\vert dt+\frac{\left(
b-x\right) ^{\mu +1}}{b-a}\int_{0}^{1}t^{\mu }\left\vert f^{\prime }\left(
tx+\left( 1-t\right) b\right) \right\vert dt.
\end{eqnarray*}%
By applying the H\"{o}lder inequality for $q>1,$ $0\leq p\leq q,$ we get%
\begin{eqnarray*}
&&\left\vert \frac{\left( x-a\right) ^{\mu }+\left( b-x\right) ^{\mu }}{b-a}%
f\left( x\right) -\frac{\Gamma \left( \mu +1\right) }{b-a}\left[
J_{x^{-}}^{\mu }f\left( a\right) +J_{x^{+}}^{\mu }f\left( b\right) \right]
\right\vert \\
&\leq &\frac{\left( x-a\right) ^{\mu +1}}{b-a}\left[ \left(
\dint\limits_{0}^{1}t^{\mu \left( \frac{q-p}{q-1}\right) }dt\right) ^{\frac{%
q-1}{q}}\left( \dint\limits_{0}^{1}t^{\mu p}\left\vert f^{\prime }\left(
tx+\left( 1-t\right) a\right) \right\vert ^{q}dt\right) ^{\frac{1}{q}}\right.
\\
&&\left. +\frac{\left( b-x\right) ^{\mu +1}}{b-a}\left(
\dint\limits_{0}^{1}t^{\mu \left( \frac{q-p}{q-1}\right) }dt\right) ^{\frac{%
q-1}{q}}\left( \dint\limits_{0}^{1}t^{\mu p}\left\vert f^{\prime }\left(
tx+\left( 1-t\right) b\right) \right\vert ^{q}dt\right) ^{\frac{1}{q}}\right]
.
\end{eqnarray*}%
It is easy to see that%
\begin{equation*}
\dint\limits_{0}^{1}t^{\mu \left( \frac{q-p}{q-1}\right) }dt=\frac{q-1}{\mu
\left( q-p\right) +q-1}.
\end{equation*}%
Hence, by $\left( \alpha ,m\right) -$logarithmically convexity of $%
\left\vert f^{\prime }\right\vert ^{q},$ we have%
\begin{eqnarray}
&&  \label{aa6} \\
&&\left\vert \frac{\left( x-a\right) ^{\mu }+\left( b-x\right) ^{\mu }}{b-a}%
f\left( x\right) -\frac{\Gamma \left( \mu +1\right) }{b-a}\left[
J_{x^{-}}^{\mu }f\left( a\right) +J_{x^{+}}^{\mu }f\left( b\right) \right]
\right\vert  \notag \\
&\leq &\frac{\left( x-a\right) ^{\mu +1}}{b-a}\left( \frac{q-1}{\mu \left(
q-p\right) +q-1}\right) ^{\frac{q-1}{q}}\left( \dint\limits_{0}^{1}t^{\mu
p}\left\vert f^{\prime }\left( x\right) \right\vert ^{t^{\alpha }}\left\vert
f^{\prime }\left( \frac{a}{m}\right) \right\vert ^{m\left( 1-t^{\alpha
}\right) }dt\right) ^{\frac{1}{q}}  \notag \\
&&+\frac{\left( b-x\right) ^{\mu +1}}{b-a}\left( \frac{q-1}{\mu \left(
q-p\right) +q-1}\right) ^{\frac{q-1}{q}}\left( \dint\limits_{0}^{1}t^{\mu
p}t^{\mu }\left\vert f^{\prime }\left( x\right) \right\vert ^{t^{\alpha
}}\left\vert f^{\prime }\left( \frac{b}{m}\right) \right\vert ^{m\left(
1-t^{\alpha }\right) }dt\right) ^{\frac{1}{q}}  \notag \\
&=&\frac{\left( x-a\right) ^{\mu +1}}{b-a}\left( \frac{q-1}{\mu \left(
q-p\right) +q-1}\right) ^{\frac{q-1}{q}}\left( \dint\limits_{0}^{1}t^{\mu
p}M^{m+t^{\alpha }\left( 1-m\right) }dt\right) ^{\frac{1}{q}}  \notag \\
&&+\frac{\left( b-x\right) ^{\mu +1}}{b-a}\left( \frac{q-1}{\mu \left(
q-p\right) +q-1}\right) ^{\frac{q-1}{q}}\left( \dint\limits_{0}^{1}t^{\mu
p}M^{m+t^{\alpha }\left( 1-m\right) }dt\right) ^{\frac{1}{q}}.  \notag
\end{eqnarray}%
If we choose $M=1,$ then 
\begin{equation*}
\dint\limits_{0}^{1}t^{\mu p}dt=\frac{1}{\mu p+1}.
\end{equation*}%
If $M<1,$ then $M^{m+t^{\alpha }\left( 1-m\right) }\leq M^{m+\alpha t\left(
1-m\right) }$, thus%
\begin{equation*}
\dint\limits_{0}^{1}t^{\mu p}M^{m+\alpha t\left( 1-m\right) }dt=\frac{%
M^{m}\left( \Gamma \left( \mu p+1\right) -\Gamma \left( \mu p+1,\ln
M^{\alpha \left( m-1\right) }\right) \right) }{\left( \ln M^{\alpha \left(
m-1\right) }\right) ^{\mu p+1}}
\end{equation*}%
Which completes the proof.
\end{proof}

\begin{corollary}
Let $f:\left[ 0,\infty \right) \rightarrow \left( 0,\infty \right) $ be \
differentiable mapping with $a,b\in \left[ 0,\infty \right) $ such that $%
a<b. $ If $\left\vert f^{\prime }\left( x\right) \right\vert ^{q}$ is $m-$%
logarithmically convex function with $\left\vert f^{\prime }\left( x\right)
\right\vert ^{q}\leq M,f^{\prime }\in L\left[ a,b\right] ,$ $m\in \left( 0,1%
\right] $ and $\mu >0,$ then the following inequality for fractional
integrals holds:%
\begin{eqnarray*}
&&\left\vert \frac{\left( x-a\right) ^{\mu }+\left( b-x\right) ^{\mu }}{b-a}%
f\left( x\right) -\frac{\Gamma \left( \mu +1\right) }{b-a}\left[
J_{x^{-}}^{\mu }f\left( a\right) +J_{x^{+}}^{\mu }f\left( b\right) \right]
\right\vert \\
&\leq &\left( \frac{q-1}{\mu \left( q-p\right) +q-1}\right) ^{\frac{q-1}{q}%
}\left( K_{2}(1,m,t)\right) ^{\frac{1}{q}}\left[ \frac{\left( x-a\right)
^{\mu +1}+\left( b-x\right) ^{\mu +1}}{\left( b-a\right) }\right]
\end{eqnarray*}%
where $q>1,$ $0\leq p\leq q$ and 
\begin{equation*}
K_{2}(1,m,t)=\left\{ 
\begin{array}{cc}
\frac{M^{m}\left( \Gamma \left( \mu p+1\right) -\Gamma \left( \mu p+1,\ln
M^{\left( m-1\right) }\right) \right) }{\left( \ln M^{\left( m-1\right)
}\right) ^{\mu p+1}} & ,M<1 \\ 
&  \\ 
\frac{1}{\mu p+1} & ,M=1%
\end{array}%
\right. .
\end{equation*}
\end{corollary}

\begin{proof}
If we set $\alpha =1$ in \ref{a5}, the proof is completed.
\end{proof}

\begin{corollary}
Let $f:\left[ 0,\infty \right) \rightarrow \left( 0,\infty \right) $ be \
differentiable mapping with $a,b\in \left[ 0,\infty \right) $ such that $%
a<b. $ If $\left\vert f^{\prime }\left( x\right) \right\vert ^{q}$ is
logarithmically convex function with $\left\vert f^{\prime }\left( x\right)
\right\vert ^{q}\leq M,f^{\prime }\in L\left[ a,b\right] $ and $\mu >0,$
then the following inequality for fractional integrals holds:%
\begin{eqnarray*}
&&\left\vert \frac{\left( x-a\right) ^{\mu }+\left( b-x\right) ^{\mu }}{b-a}%
f\left( x\right) -\frac{\Gamma \left( \mu +1\right) }{b-a}\left[
J_{x^{-}}^{\mu }f\left( a\right) +J_{x^{+}}^{\mu }f\left( b\right) \right]
\right\vert \\
&=&\left( \frac{q-1}{\mu \left( q-p\right) +q-1}\right) ^{\frac{q-1}{q}%
}\left( \frac{1}{\mu p+1}\right) ^{\frac{1}{q}}\left[ \frac{\left(
x-a\right) ^{2}+\left( b-x\right) ^{2}}{\left( b-a\right) }\right]
\end{eqnarray*}%
where $q>1,$ $0\leq p\leq q$.
\end{corollary}

\begin{proof}
If we set $\alpha =m=1$ in \ref{aa6}, the proof is completed.
\end{proof}

\begin{corollary}
Let $f:\left[ 0,\infty \right) \rightarrow \left( 0,\infty \right) $ be \
differentiable mapping with $a,b\in \left[ 0,\infty \right) $ such that $%
a<b. $ If $\left\vert f^{\prime }\left( x\right) \right\vert ^{q}$ is $%
\left( \alpha ,m\right) -$logarithmically convex function with $\left\vert
f^{\prime }\left( x\right) \right\vert ^{q}\leq M,f^{\prime }\in L\left[ a,b%
\right] $ and $\left( \alpha ,m\right) \in \left( 0,1\right] \times \left(
0,1\right] ,$ then the following inequality holds:%
\begin{eqnarray}
&&\left\vert f\left( x\right) -\frac{1}{b-a}\int_{a}^{b}f(u)du\right\vert
\label{a6} \\
&\leq &\left( \frac{q-1}{2q-p-1}\right) ^{\frac{q-1}{q}}\left( K_{1}(\alpha
,m,t)\right) ^{\frac{1}{q}}\left[ \frac{\left( x-a\right) ^{2}+\left(
b-x\right) ^{2}}{\left( b-a\right) }\right]  \notag
\end{eqnarray}%
where $q>1,$ $0\leq p\leq q$ and 
\begin{equation*}
K_{3}(\alpha ,m,t)=\left\{ 
\begin{array}{cc}
\frac{M^{m}\left( \Gamma \left( p+1\right) -\Gamma \left( p+1,\ln M^{\alpha
\left( m-1\right) }\right) \right) }{\left( \ln M^{\alpha \left( m-1\right)
}\right) ^{p+1}} & ,M<1 \\ 
&  \\ 
\frac{1}{p+1} & ,M=1%
\end{array}%
\right. .
\end{equation*}
\end{corollary}

\begin{proof}
If we set $\mu =1$ in \ref{a5}, the proof is completed.
\end{proof}

\begin{corollary}
Let $f:\left[ 0,\infty \right) \rightarrow \left( 0,\infty \right) $ be \
differentiable mapping with $a,b\in \left[ 0,\infty \right) $ such that $%
a<b. $ If $\left\vert f^{\prime }\left( x\right) \right\vert ^{q}$ is $%
\left( \alpha ,m\right) -$logarithmically convex function with $\left\vert
f^{\prime }\left( x\right) \right\vert ^{q}\leq M,f^{\prime }\in L\left[ a,b%
\right] $ and $\left( \alpha ,m\right) \in \left( 0,1\right] \times \left(
0,1\right] ,$ then the following inequality holds:%
\begin{eqnarray*}
&&\left\vert f\left( x\right) -\frac{1}{b-a}\int_{a}^{b}f(u)du\right\vert \\
&\leq &\left( \frac{1}{2}\right) ^{\frac{q-1}{q}}\left( K_{1}(\alpha
,m,t)\right) ^{\frac{1}{q}}\left[ \frac{\left( x-a\right) ^{2}+\left(
b-x\right) ^{2}}{\left( b-a\right) }\right]
\end{eqnarray*}%
where $q>1,$ $0\leq p\leq q$ and 
\begin{equation*}
K_{4}(\alpha ,m,t)=\left\{ 
\begin{array}{cc}
\frac{M^{m}\left( \Gamma \left( 2\right) -\Gamma \left( 2,\ln M^{\alpha
\left( m-1\right) }\right) \right) }{\left( \ln M^{\alpha \left( m-1\right)
}\right) ^{2}} & ,M<1 \\ 
&  \\ 
\frac{1}{2} & ,M=1%
\end{array}%
\right. .
\end{equation*}
\end{corollary}

\begin{proof}
If we set $p=1$ in \ref{a6}, the proof is completed.
\end{proof}


\begin{thebibliography}{99}
\bibitem{1} R.-F. Bai, F. Qi and B.-Y. Xi, \textit{Hermite-Hadamard type
inequalities for the }$m-$\textit{\ and }$\left( \alpha ,m\right) -$\textit{%
logarithmically convex functions}, Filomat 27 (2013), 1-7.

\bibitem{BPR} M.K. Bakula, J. Pe\v{c}ari\'{c} and M. Ribibi\'{c}, \textit{%
Companion inequalities to Jensen's inequality for }$m-$\textit{convex and }$%
(\alpha ,m)-$\textit{convex functions}, J. Inequal. Pure and Appl. Math., 7
(5) (2006), Article 194.

\bibitem{ST} S.S. Dragomir and G. Toader, \textit{Some inequalities for }$m-$%
\textit{convex functions}, Studia University Babes Bolyai, Mathematica, 38
(1) (1993), 21-28.

\bibitem{MIH} V.G. Mihe\c{s}an, \textit{A generalization of the convexity},
Seminar of Functional Equations, Approx. and Convex, Cluj-Napoca (Romania)
(1993).

\bibitem{TOA} G. Toader,\textit{\ Some generalization of the convexity},
Proc. Colloq. Approx.\textit{\ Opt.}, Cluj-Napoca, (1984), 329-338.

\bibitem{TOA2} G. Toader, \textit{On a generalization of the convexity},
Mathematica, 30 (53) (1988), 83-87.

\bibitem{10} A. Ostrowski,\textit{\ \"{U}ber die Absolutabweichung einer
differentierbaren Funktion von ihren Integralmittelwert}, Comment. Math.
Helv., 10, 226-227, (1938).

\bibitem{a1} M.E. \"{O}zdemir, H. Kavurmaci, E. Set, \textit{Ostrowski's
type inequalities for }$(\alpha ,m)-$\textit{convex functions}, KYUNGPOOK
Math. J. 50 (2010) 371--378.

\bibitem{2} H. Kavurmaci, M. Avci and M.E. \"{O}zdemir, \textit{New
Ostrowski type inequalities for }$m-$\textit{convex functions and
applications}, Hacettepe Journal of Mathematics and Statistics, Volume 40
(2) (2011), 135 -- 145.

\bibitem{11} M. Alomari and M. Darus, \textit{Some Ostrowski type
inequalities for convex functions with applications}, RGMIA Res. Rep. Coll.,
(2010) 13, 2, Article 3. [ONLINE: http://ajmaa.org/RGMIA/v13n2.php].

\bibitem{a2} E. Set, \textit{New inequalities of Ostrowski type for mappings
whose derivatives are }$s$\textit{-convex in the second sense via fractional
integrals}, Comput. Math. Appl., 63 (2012) 1147-1154.

\bibitem{anastas} S. Belarbi and Z. Dahmani, \textit{On some new fractional
integral inequalities}, J. Ineq. Pure and Appl. Math., 10(3), Art. 86 (2009).

\bibitem{dahmani} Z. Dahmani, \textit{New inequalities in fractional
integrals,} International Journal of Nonlinear Science, 9(4), 493-497 (2010).

\bibitem{zdahm} Z. Dahmani, \textit{On Minkowski and Hermite-Hadamard
integral inequalities via fractional integration}, Ann. Funct. Anal. 1(1),
51-58 (2010).

\bibitem{zdah} Z. Dahmani, L. Tabharit and S. Taf, \textit{Some fractional
integral inequalities}, Nonl. Sci. Lett. A., 1(2), 155-160 (2010).

\bibitem{zeki2} M.Z. Sar\i kaya, E. Set, H. Yaldiz and N. Ba\c{s}ak, \textit{%
Hermite-Hadamard's inequalities for fractional integrals and related
fractional inequalities}, Mathematical and Computer Modelling, In Press.

\bibitem{dahtab} Z. Dahmani, L. Tabharit and S. Taf, \textit{New
generalizations of Gr\"{u}ss inequality using Riemann-Liouville fractional
integrals}, Bull. Math. Anal. Appl., 2(3), 93-99 (2010).

\bibitem{a} M.E. \"{O}zdemir, H. Kavurmac\i\ and M. Avc\i , \textit{New
inequalities of Ostrowski type for mappings whose derivatives are }$(\alpha
,m)$\textit{-convex via fractional integrals}, RGMIA Research Report
Collection, 15, Article 10, 8 pp (2012).
\end{thebibliography}
\end{document}